\documentclass[graybox]{svmult}

\usepackage{mathptmx}        
\usepackage{makeidx}         
\usepackage{cite}
\usepackage{graphicx}        
\usepackage{grffile}
\usepackage{epstopdf}        
\usepackage{amsmath}
\usepackage{multicol}        
\usepackage[bottom]{footmisc}

\makeindex             

\begin{document}


\title*{New families of periodic orbits for the planar three-body problem computed with high precision}

\author{I. Hristov$^{1,a}$, R. Hristova$^1$, I. Puzynin$^2$,  T. Puzynina$^2$, Z. Sharipov$^{2,b}$, Z. Tukhliev$^2$\\
\vspace{0.5cm}
\emph{$^1$ Faculty of Mathematics and Informatics, Sofia University, Sofia, Bulgaria}\\
\emph{$^2$ Meshcheryakov Laboratory of Information Technologies, JINR,  Dubna, Russia}\\
\vspace{0.5cm}
\emph{E-mails:  $^a$  ivanh@fmi.uni-sofia.bg \hspace{0.3cm} $^b$ zarif@jinr.ru}}
\authorrunning{I. Hristov, R. Hristova, I. Puzynin et al.}
\titlerunning{New families of periodic orbits for the planar three-body problem }
\maketitle

\abstract {In this paper we use a Modified Newton's method based on the Continuous analog of Newton's method and high precision arithmetic for a general numerical search of periodic orbits for the planar three-body problem. We consider relatively short periods and a relatively coarse search-grid. As a result, we found 123 periodic solutions belonging to 105 new topological families that are not included in the database in [Science China Physics, Mechanics \& Astronomy 60.12 (2017)]. The extensive numerical search is achieved by a parallel solving of many independent tasks using many cores in a computational cluster.}

\vspace{0.5 cm}


\section{Introduction}

Numerical search of periodic orbits is an important tool for further investigation and understanding of the three-body problem. In 2013 Shuvakov and Dmitrashinovich found 13 new topological families applying a numerical algorithm based on the gradient descent method in the standard double precision arithmetic  \cite{Shuv1, Shuv2}.
It is well known that the three-body problem is very sensitive on the initial conditions. This means that working with double precision strongly limits the number of solutions that can be found. This limitation was recognized by Li and Liao, in 2017 they applied Newton's method with high precision arithmetic to find more than 600 new families of periodic orbits  \cite{Liao1}. In this paper we use a Modified Newton's method based on the Continuous analog of Newton's method and high precision arithmetic for a general numerical search of periodic orbits for the planar three-body problem. By considering relatively short periods and a relatively coarse search-grid, we found 123 solutions belonging to 105 new topological families that are not included in \cite{Liao1}.

\section{Mathematical model}
We consider three bodies with equal masses. They are treated as mass points.
A planar motion of the three bodies is considered.
The normalized differential equations describing the motion of the bodies are:
\begin{equation}
\ddot{r}_i=\sum_{j=1,j\neq i}^{3} \frac{(r_j-r_i)}{{\|r_i -r_j\|}^3}, i=1,2,3.
\end{equation}
The vectors $r_i$, $\dot{r}_i$ have two components: $r_i=(x_i, y_i), \dot{r}_i=(\dot{x}_i, \dot{y}_i).$
The system (1) can be written as a first order one this way:
\begin{equation}
\dot{x}_i={vx}_i, \hspace{0.1 cm} \dot{y}_i={vy}_i,  \hspace{0.1 cm} \dot{vx}_i=\sum_{j=1,j\neq i}^{3} \frac{(x_j-x_i)}{{\|r_i -r_j\|}^3},  \hspace{0.1 cm} \dot{vy}_i=\sum_{j=1,j\neq i}^{3} \frac{(y_j-y_i)}{{\|r_i -r_j\|}^3}, \hspace{0.1 cm} i=1,2,3
\end{equation}
We numerically solve the problem in this first order form. Hence, we have a vector of 12 unknown functions
$ u(t)={(x_1, y_1, {vx}_1, {vy}_1, x_2, y_2, {vx}_2, {vy}_2, x_3, y_3, {vx}_3, {vy}_3)}^\top$.
We search for periodic planar collisionless orbits as in \cite{Shuv1, Shuv2, Liao1}: with zero angular momentum and symmetric initial
configuration with parallel velocities:
\begin{equation}
\begin{aligned}
(x_1(0),y_1(0))=(-1,0), \hspace{0.2 cm} (x_2(0),y_2(0))=(1,0), \hspace{0.2 cm} (x_3(0),y_3(0))=(0,0) \\
({vx}_1(0),{vy}_1(0))=({vx}_2(0),{vy}_2(0))=(v_x,v_y) \hspace{2 cm}\\
({vx}_3(0),{vy}_3(0))=-2({vx}_1(0),{vy}_1(0))=(-2v_x, -2v_y) \hspace{1.5 cm}
\end{aligned}
\end{equation}
Here the velocities $v_x\in [0,1], v_y\in [0,1]$ are parameters. The periods of the orbits are denoted with $T$. So, our goal is to find triplets $(v_x, v_y, T)$
for which the periodicity condition $u(T)=u(0)$ is fulfilled.

\section{Description of the numerical procedure}
The numerical procedure consists of three stages. During stage (I) we compute initial approximations for the correction method (Modified Newton's method) by applying the grid-search method on the rectangle  $[0,0.8]\times[0,0.8]$ for $v_x, v_y$ with a stepsize $1/2048$.
The candidates for correction are the triplets $(v_x, v_y, T)$, such that the return proximity $P(T)$ has local minima
on the grid for $v_x, v_y$ and $P(T)$ is less than 0.7: $P(T) = \min_{1<t\leq T_0}P(t) < 0.7.$ Here $P(t)$ is defined  this way:
$$P(t)=\sqrt{\sum_{i=1}^{3}{\|r_i(t)-r_i(0)\|}^{2}+\sum_{i=1}^{3}{\|\dot{r}_i(t)-\dot{r}_i(0)\|}^{2}}$$

For every grid point the system (2) with initial conditions (3) is solved numerically up to  $T_0=70$ by using multiple precision Taylor series method (MPTSM)  \cite{Barrio} with order 154 and precision 134 decimal digits. During stage (II) we apply the modified Newton's method, which can converge or diverge. Convergence means that we catch a periodic solution. At each iteration step of the modified Newton's method, a system of 36 ODEs (the original differential equations plus the differential equations for the partial derivatives with respect to the parameters $v_x$ and $v_y$ ) is solved. MPTSM is used again with the same order - 154 and the same precision - 134 decimal digits. During stage (III) we apply the classic Newton's method in order to specify the solutions with more correct digits (150 correct digits in this work) and verify them. We use MPTSM with increased order and precision for solving the system of 36 ODEs. During this stage we make two computations. The first computation is with 220-th order method and 192 decimal digits of precision and the second computation is for verification - with 264-th order method and 231 digits of precision. All the details of the Newton's and the modified Newton's method for the planar three-body problem can be seen in our recent work \cite{Ourpaper}.
The GMP library (the GNU Multiple Precision arithmetic library) \cite{GMP} is used for floating point arithmetic in our C-programs.

The extensive computations are performed in "Nestum" cluster, Sofia, Bulgaria \cite{nestum} and "Govorun" supercomputer, Dubna, Russia \cite{HybriLIT}. The two platforms use SLURM as a cluster management and job scheduling system (batch system).  The input data file for stage (I) consists of the set of all grid-points $(v_x, v_y)$ in the search window $[0,0.8]\times[0,0.8]$ with a stepsize $1/2048$. The input data file for stage (II) consists of the set of all triplets $(v_x, v_y, T)$, obtained from the first stage (all candidates for correction with the modified Newton's method). At last, the input data file for the third stage consists of the set of triplets $(v_x, v_y, T)$ for which the modified Newton's method from stage (II) converges (all caught solutions).  Each grid-point $(v_x, v_y)$ at stage (I) or each triplet at stage (II) and (III) can be processed independently in parallel. In fact, we have an example of the so-called embarrassing parallelism where no communication between tasks is needed. The computational process consists of dividing the input data file into many files and submitting the corresponding jobs associated with each data file to the batch system. When the jobs are done, the output files are assembled into one output file, which is eventually additionally processed. Shell scripts are used to automatize the process of file distribution, the jobs submission and to gather the results.  At all three stages of the numerical procedure these scripts are applied.  We achieve a substantial speedup, which is absolutely needed in order to solve the problem in a foreseeable time.

\section{Results}
For each found solution we compute the free group element (its topological family) in order to classify it \cite{Shuv1,Shuv2}.
Also we give the four numbers $(v_x,v_y,T,T^*)$ with 150 correct digits,
where $T^*$ is the scale-invariant period. The scale-invariant period $T^*$ is defined as  $T^{*}=T{|E|}^{\frac{3}{2}},$
where  E is the energy of our initial configuration: $E=-2.5 + 3({v_x}^2+{v_y}^2)$. Equal $T^*$ for two different initial conditions means two different representations of the same solution.
We found 123 solutions belonging to 105 new topological families that are not included in \cite{Liao1} and also are not their satellites. The initial velocities of the 695 solutions from \cite{Liao1} (blue points) and the initial velocities of our 123 solutions (red points) are shown in Figure 1. As it is seen from this figure, our solutions (red points) are found in regions where old solutions (blue points) are relatively sparse. The four numbers $(v_x, v_y, T, T^*)$ with 150 correct digits, the free group elements, their lengths and animations in the real space for these 123 solutions can be found in \cite{rada3body}.
Our numerical procedure is also able to find all $124$ solutions with periods $T < 70$ from \cite{Liao1}, except for two of them. Let us mention that a finer search-grid is used in \cite{Liao1} (with a stepsize $1/4000$) and periods there are up to 200. The results show that the modified Newton's method based on the Continuous analog of Newton's method has a great potential for search of new periodic orbits of the three-body problem. We did additional numerical tests for the stability of the presented solutions. These tests (although not fully rigorous) show that all new-found solutions are most possibly unstable.
\begin{figure}
\begin{center}
\includegraphics[scale=0.5]{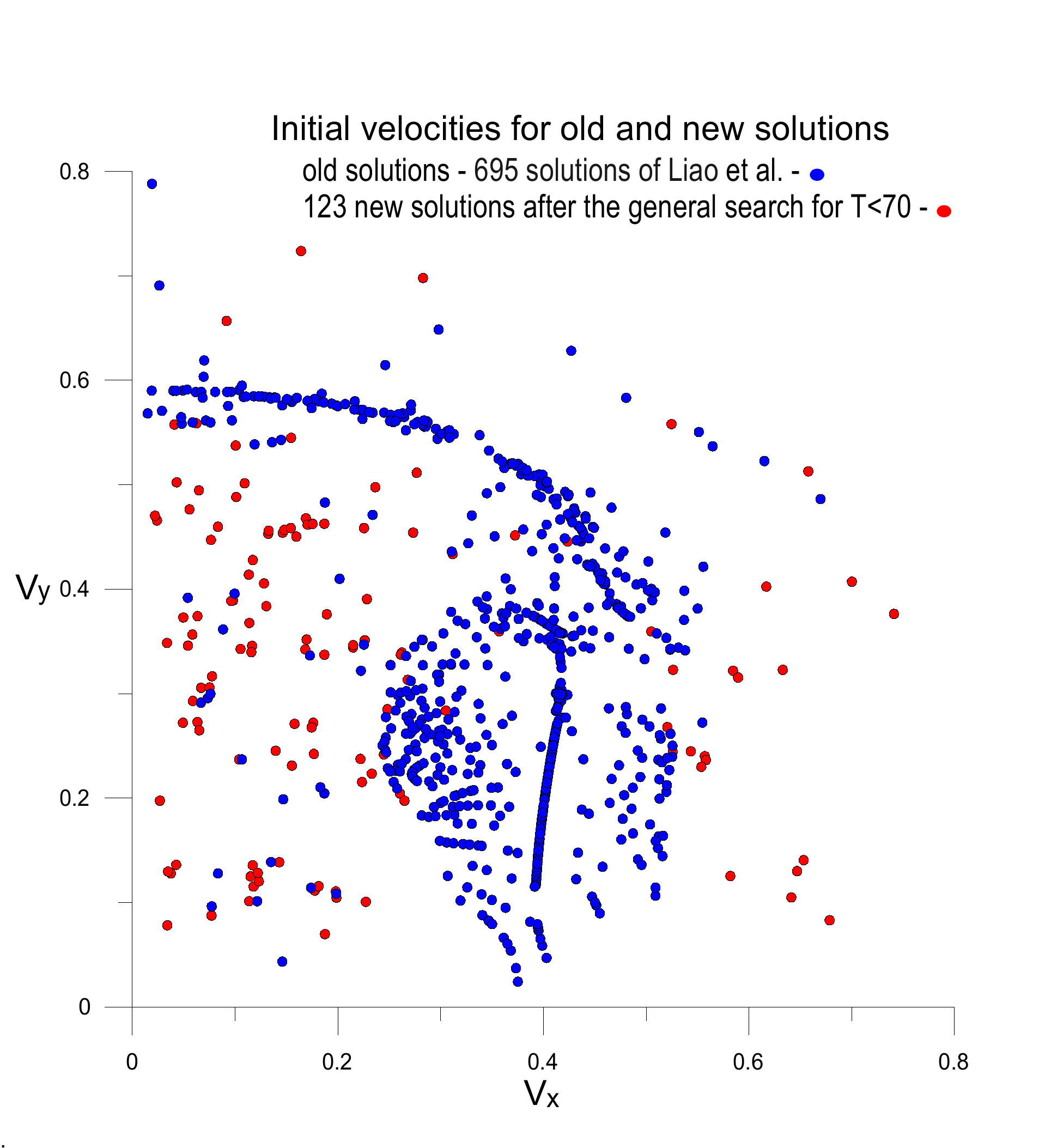}
\caption{Initial velocities for all new-found  123 solutions (red points) and 695 old solutions of Liao et al. (blue points)}
\label{fig:1}
\end{center}
\end{figure}

\begin{acknowledgement}

We thank for the opportunity to use the computational resources of the Nestum cluster, Sofia, Bulgaria
and  the "Govorun" supercomputer at the Meshcheryakov Laboratory of Information Technologies of JINR, Dubna, Russia.
The work is supported by a grant of the Plenipotentiary Representative of the Republic of Bulgaria at JINR, Dubna, Russia.

\end{acknowledgement}

\end{document}